\documentclass[a4,12pt,fleqn]{article}
\usepackage{amsmath,amsthm,enumerate,amssymb}

\begin{document}

\newcommand{\ep}{\varepsilon}

\renewcommand{\thefootnote}{\fnsymbol{footnote}}

\begin{center}
\begin{Large}
{\bf 
A fundamental inequality between products of $x^p - 1$. 
}

\vskip1cm 

\textsc{Keiichi Watanabe}
\end{Large}
\end{center} 

\vskip1cm 

\textsc{
Keiichi Watanabe: 
Department of Mathematics, Faculty of Science, Niigata University, 
950-2181, Japan}

\textit{E-mail address}: wtnbk@math.sc.niigata-u.ac.jp 

\vskip1cm

\noindent 
{\bf Abstract.} 
We will show certain functional inequalities between some products of $x^p - 1$.



\vskip1cm

\section{Introduction}

\bigskip 

This article is concerned with an elementary topic on elementary functions. 

It is easy to see the inequalities 
\begin{align*}
595 \left(x^{6} - 1\right)\left(x^{8} - 1\right)^2\left(x^{9} - 1\right) 
\le 
1728 \left(x^{2} - 1\right)\left(x^{5} - 1\right)\left(x^{7} - 1\right)\left(x^{17} - 1\right)
\end{align*}
or 
\begin{align*}
48 \left(x^2 - 1\right)\left(x^3 - 1\right)\left(x^5 - 1\right)\left(x^7 - 1\right)\left(x^{11} - 1\right) 
\le 
385 \left(x - 1\right)^2\left(x^4 - 1\right)^2\left(x^{18} - 1\right) 
\end{align*}
for arbitrary $1 < x$, 
if they are provided as the matter to be proved. 
However, if we would like to estimate functions of the form 
\begin{align*} 
\prod\left(x^{p_j} - 1\right) 
\end{align*} 
by simpler ones, how can we guess what forms and coefficients are possible? 

\bigskip 
\noindent 
{\bf Example.} The following inequality does not hold on an interval contained in $1 < x$. 
\begin{align*} 
225(x^{2} - 1)^2(x^{8} - 1)^2 \le 256(x - 1)(x^{5} - 1)^2(x^{9} - 1). 
\end{align*} 

In this paper, we will prove certain functional inequalities 
which might be related to the problem as mentioned above, 
although the efficiency and possible applications to other branches of mathematics are still to be clarified.  

The following theorem is a starting point to the main result in this paper, 
whose proof is based on an operator inequality by Furuta [1] 
and the argument related to the best possibility of that by Tanahashi [2].  

\bigskip 
\noindent 
{\bf Theorem 1.1 ([3]).} 
Let $0 \le p, \, 1 \le q$ and $0 \le r$ with $p + r \le (1 + r)q$. 
If $0 < x$, then 
\begin{align*}
x^{\frac{1 + r - \frac{p + r}{q}}{2}}\left(x^p - 1\right)\left(x^{\frac{p + r}{q}} - 1\right) 
\le \, \frac{p}{q} \left(x^{p + r} - 1\right)\left(x - 1\right). 
\end{align*}

\bigskip 
\noindent
{\bf Proposition 1.2.} 
Let $1 \le p, \, 0 \le r$. 
Then, for arbitrary $0 < x$, 
\begin{align}
(p + r)(x^{p} - 1)(x^{1 + r} - 1) 
\le 
p(1 + r)(x - 1)(x^{p + r} - 1). 
\end{align} 

\bigskip 
\noindent 
Proof. 
Put $\displaystyle q = \frac{p + r}{1 + r}$. Since $1 \le p$, we have $1 \le q$ and hence Proposition 1.2 

\smallskip
\noindent
immediately follows from Theorem 1.1. 
$\qed$ 

\bigskip 

Although it is likely there exists an elementary and direct proof of Proposition 1.2 
without using matrix inequality, the author has not found it yet.  

\bigskip 
\noindent
{\bf Theorem 1.3.} 
Let $0 < p_1 \le p_2, \enskip 0 < q_1 \le q_2, 
\enskip p_1 + p_2 = q_1 + q_2$ and $p_2 \le q_2$. 
Then, for arbitrary $0 < x$, 
\begin{align}
q_1q_2(x^{p_1} - 1)(x^{p_2} - 1) 
\le 
p_1p_2(x^{q_1} - 1)(x^{q_2} - 1) 
\end{align} 

\bigskip 
\noindent 
Proof. 
For a moment, we add $1 \le p_1, \, p_2, \, q_2$ and $q_1 = 1$ to the assumption. 
Apply Proposition 1.2 with $p = p_1, \, r = p_2 - 1$, 
then the inequality (1) implies 
\begin{align*}
q_2(x^{p_1} - 1)(x^{p_2} - 1) 
\le 
p_1p_2(x - 1)(x^{q_2} - 1). 
\end{align*} 

In general, note that $q_1 \le p_1$. 
Dividing by $q_1$, we have $\displaystyle 1 \le \frac{p_1}{q_1} \le \frac{p_2}{q_1}, \enskip 1 \le \frac{q_2}{q_1}$, 
\begin{align*} 
\frac{p_1}{q_1} + \frac{p_2}{q_1} = 1 + \frac{q_2}{q_1} 
\quad \mbox{and} \quad 
\frac{p_2}{q_1} \le \frac{q_2}{q_1}. 
\end{align*}  

By the first part of the proof, 
\begin{align*}
&\frac{q_2}{q_1} (x^{\frac{p_1}{q_1}} - 1)(x^{\frac{p_2}{q_1}} - 1) 
\le 
\frac{p_1}{q_1} \cdot \frac{p_2}{q_1} (x - 1)(x^{\frac{q_2}{q_1}} - 1) 
\end{align*} 
for arbitrary $0 < x$. 
By substituting $x^{q_1}$ to $x$ in the above inequality, 
it is immediate to see the inequality (2). 
$\qed$

\section{Main Result}

\bigskip
\noindent
{\bf Theorem 2.1.} 
Let $n$ be a natural number. 
Suppose $0 < p_1 \le \cdots \le p_n, \enskip 0 < q_1 \le \cdots \le q_n$ and 
\begin{align*} 
\sum_{j=1}^n p_j &= \sum_{j=1}^n q_j \\[+5pt] 
\sum_{j=2}^n p_j &\le \sum_{j=2}^n q_j \\ 
&\vdots \\ 
p_{n - 1} + p_n &\le q_{n - 1} + q_n \\[+5pt] 
p_n &\le q_n. 
\end{align*} 
Then, for arbitrary $1 < x$, 
\begin{align}
&\prod_{j=1}^n q_j \prod_{j = 1}^{n}(x^{p_j} - 1) 
\le 
\prod_{j=1}^n p_j \prod_{j = 1}^{n}(x^{q_j} - 1). 
\end{align} 
If $n$ is even, the inequality (3) holds for arbitrary $0 < x < 1$. 
If $n$ is odd, the reverse inequality of (3) holds for arbitrary $0 < x < 1$.

\bigskip 
\noindent
Proof. 
The step $2$ is exactly Theorem 1.3 in the previous section. 
Assume that the step $n$ is valid. 
Suppose $0 < p_1 \le \cdots \le p_n \le p_{n + 1}, \enskip 0 < q_1 \le \cdots \le q_n \le q_{n + 1}$ and 
\begin{align*} 
\sum_{j=1}^{n + 1} p_j &= \sum_{j=1}^{n + 1} q_j \\[+5pt] 
\sum_{j=2}^{n + 1} p_j &\le \sum_{j=2}^{n + 1} q_j \\ 
&\vdots \\ 
p_n + p_{n + 1} &\le q_n + q_{n + 1} \\[+5pt] 
p_{n + 1} &\le q_{n + 1}. 
\end{align*}  

\bigskip 

There exists a number $k$ such that $1 \le k \le n$ and 
\begin{align*} 
q_1 \le \cdots \le q_k \le p_1 \le q_{k + 1} \le \dots \le q_{n + 1}. 
\end{align*} 
Take a real number $q^\prime$ which is determined by 
$q_k + q_{k + 1} = p_1 + q^\prime$. Then 
\[ 
q_k \le q^\prime = q_k + q_{k + 1} - p_1 \le q_{k + 1}. 
\] 

By the step $2$, 
\begin{align}
q_kq_{k + 1}(x^{p_1} - 1)(x^{q^\prime} - 1) \le p_1q^\prime(x^{q_k} - 1)(x^{q_{k + 1}} - 1). 
\end{align} 

Since 
\begin{align*} 
p_1 + \cdots + p_{n + 1} = q_1 + \cdots + q_{n + 1} = p_1 + q^\prime + \sum_{j \ne k, k + 1}q_j, 
\end{align*} 
we have 
\begin{align*} 
p_2 + \cdots + p_{n + 1} = q^\prime + \sum_{j \ne k, k + 1}q_j 
\end{align*} 
and 
\[ 
q_1 \le \cdots \le q_{k - 1} \le q^\prime \le q_{k + 2} \le \cdots \le q_{n + 1}. 
\] 

The previous equality and $q_1 \le p_2$ yield 
\begin{align*} 
p_3 + \cdots + p_{n + 1} \le q^\prime + \sum_{j \ge 2, j \ne k, k + 1}q_j. 
\end{align*} 

It follows from the assumption that 
\begin{align*} 
p_3 + p_4 + \cdots + p_{n + 1} 
&\le q_3 + q_4 + \cdots + q_{n + 1} \\[+5pt] 
&= p_1 + q^\prime + \sum_{j \ge 3, j \ne k, k + 1}q_j 
\le p_3 + q^\prime + \sum_{j \ge 3, j \ne k, k + 1}q_j, 
\end{align*}  
and hence 
\begin{align*} 
p_4 + \cdots + p_{n + 1} 
\le q^\prime + \sum_{j \ge 3, j \ne k, k + 1}q_j. 
\end{align*} 

Similarly we have 
\begin{align*} 
p_5 + \cdots + p_{n + 1} 
&\le q^\prime + \sum_{j \ge 4, j \ne k, k + 1}q_j \\ 
&\vdots \\ 
p_{k + 1} + \cdots + p_{n + 1} &\le q^\prime + \sum_{j = k + 2}^{n + 1} q_j 
\end{align*} 
\begin{align*} 
p_{k + 2} + \cdots + p_{n + 1} &\le q_{k + 2} + \cdots + q_{n + 1} \\ 
&\vdots \\ 
p_n + p_{n + 1} &\le q_n + q_{n + 1} \\ 
p_{n + 1} &\le q_{n + 1}. 
\end{align*}

Therefore, $n$-tuples $\{p_2, \cdots, p_{n + 1}\}, \, \{q_1, \cdots, q_{k - 1}, q^\prime, q_{k + 2}, \cdots, q_{n + 1}\}$ 
satisfy the assumption of the step $n$, and so we can obtain 
\begin{align}
&q^\prime \prod_{j \ne k, k + 1} q_j \prod_{j = 2}^{n + 1}(x^{p_j} - 1) 
\le 
\prod_{j=2}^{n + 1} p_j (x^{q^\prime} - 1) \prod_{j \ne k, k + 1}(x^{q_j} - 1) 
\end{align}  
for arbitrary $1 < x$. 

From (4) and (5), it is immediate to see that 
\begin{align*}
&\prod_{j=1}^{n + 1} q_j \prod_{j = 1}^{n + 1}(x^{p_j} - 1) 
\le 
\prod_{j=1}^{n + 1} p_j (x - 1) \prod_{j = 1}^{n + 1}(x^{q_j} - 1) 
\end{align*} 
for $1 < x$. 

The last assertion of the theorem can be easily seen 
by substituting $\displaystyle \frac1{x}$ for $0 < x < 1$ and 
multiplying $\displaystyle x^{p_1 + \cdots + p_n}$ to both sides. 

This completes the proof. 
$\qed$

\bigskip 
\bigskip 
\noindent 
{\bf Remark 2.2.} Each following example of the case $n = 5$ 
does not satisfy one of the conditions for parameters in the assumption of Theorem 2.1 
and the inequality (3) does not hold for all $1 < x$.  
\begin{enumerate} 
\item[(1)] $p_5 > q_5$. 
\begin{align*} 
4 \cdot 6 \cdot 8(x^{2} - 1)^2(x^{3} - 1)^2(x^{10} - 1) 
\le 2^2 \cdot 3^2 \cdot 10(x - 1)^2(x^{4} - 1)(x^{6} - 1)(x^{8} - 1) 
\end{align*} 
\item[(2)] $p_4 + p_5 > q_4 + q_5$. 
\begin{align*} 
6 \cdot 8 \cdot 12(x^{2} - 1)^3(x^{11} - 1)^2
\le 2^3 \cdot 11^2(x - 1)^2(x^{6} - 1)(x^{8} - 1)(x^{12} - 1) 
\end{align*} 
\item[(3)] $p_3 + p_4 + p_5 > q_3 + q_4 + q_5$. 
\begin{align*} 
5^2 \cdot 10^2(x^{2} - 1)^2(x^{9} - 1)^3 
\le 2^2 \cdot 9^3(x - 1)(x^{5} - 1)^2(x^{10} - 1)^2 
\end{align*} 
\item[(4)] $p_2 + p_3 + p_4 + p_5 > q_2 + q_3 + q_4 + q_5$. 
\begin{align*} 
3^4 \cdot 9(x - 1)(x^{5} - 1)^4 
\le 5^4(x^{3} - 1)^4(x^{9} - 1). 
\end{align*} 
\end{enumerate}

\bigskip 
\noindent 
{\bf Remark 2.3.} 
There exists an example of the case $n = 3$ such that 
$p_3 > q_3$ but the inequality (3) holds for $1 < x$. 
\begin{align*} 
5 \cdot 6(x^{2} - 1)(x^{3} - 1)(x^{7} - 1) \le 2 \cdot 3 \cdot 7(x - 1)(x^{5} - 1)(x^{6} - 1) 
\end{align*} 
Indeed, 
\begin{align*} 
&2 \cdot 3 \cdot 7(x - 1)(x^{5} - 1)(x^{6} - 1) - 5 \cdot 6(x^{2} - 1)(x^{3} - 1)(x^{7} - 1) \\[+5pt]
&= 6 (x - 1)^7 (x + 1) (x^2 + x + 1) (2 x^2 + 3 x + 2). 
\end{align*}

\section*{Acknowledgements}
The author was supported in part by 
Grants-in-Aid for Scientific Research, 
Japan Society for the Promotion of Science.

\renewcommand{\refname}{References}

\end{document}